\documentclass[11pt]{amsart}
\usepackage{amsmath, amssymb, amsthm}
\usepackage[margin=1in]{geometry}
\usepackage{mathrsfs}
\usepackage{xcolor}
\usepackage{hyperref}
\usepackage{dsfont}
\usepackage{comment}

\DeclareMathOperator{\supp}{supp}

\newcommand{\R}{\mathbb{R}}
\newcommand{\N}{\mathbb{N}}

\newtheorem{theorem}{Theorem}
\newtheorem{definition}{Definition}
\newtheorem{remark}{Remark}

\title[A generalized Liouville theorem]{A generalized Liouville theorem via division}

\author{
  David Lee
}
\address{School of Mathematical and Physical Sciences, University of Technology Sydney (UTS)}
\email{davidchanwoo.lee@uts.edu.au}
\date{\today}

\begin{document}

\begin{abstract}
We study the equation $P(i\nabla)u=0$ on $\R^d$ for a class of admissible symbols $P$ whose zero set is the unit sphere $S^{d-1}$ and which
vanish there to some finite order. Working in the framework of Lizorkin
distributions, and hence without any boundedness or decay hypothesis on $u$, we give a
complete classification of the solutions: $u$ solves $P(i\nabla)u=0$ if and only if
$\hat{u}$ is a multi-layer distribution on $S^{d-1}$ of order at most $N$. Alternatively, $u$ solves $P(i\nabla)u=0$ if and only if $(1+\Delta)^{N+1}u=0$ if $P$ satisfies a flatness condition. The proof
recasts the equation as a division problem and combines the order of vanishing of $P$
with the structure theorem for distributions. This unifies and
extends known Helmholtz-type rigidity results, which correspond to a simple zero on the
sphere, to symbols with zeros of arbitrary finite order. 
\end{abstract}

\subjclass[2010]{35A08, 35JXX, 35R11, 42B37, 42B10,42B20}
\thanks{Keywords: harmonic analysis; Liouville theorem; non-local}
\maketitle

\section{Introduction}

The study of Helmholtz-type equations has, in recent years, been shaped by a striking
rigidity phenomenon: the solutions of a nonlocal equation are often forced to coincide
with those of its classical, local counterpart. This line of inquiry was initiated by
Fall and Weth~\cite{MR3511811} and Guan, Murugan and Wei~\cite{MR4546886}. Specifically,
for $0<s<1$, Guan, Murugan and Wei showed that the bounded solutions of the fractional
Helmholtz equation $(-\Delta)^s u = u$ in $\R^d$ that decay to $0$ are precisely the
bounded solutions of the classical equation $-\Delta u = u$; the one-dimensional case (for bounded and non-decaying solutions)
had earlier been treated by Fall and Weth~\cite{MR3511811}.

This rigidity has since been extended in several
directions~\cite{MR4709301, MR4845366, hauer2023characterization, MR3477075, Lee2026}.
The most general such result is due to Cheng, Li and Yang~\cite{MR4709301}, who proved
that any $u\in L^\infty(\R^d)$ solving
\begin{equation*}
	\Phi(-\Delta)u = \Phi(1)u, \quad \text{in $\mathcal{S}'(\R^d)$,}
\end{equation*}
necessarily satisfies $-\Delta u = u$, for a broad class of symbols $\Phi$ that
encompasses both Bernstein functions of the Laplacian and the logarithmic Laplacian (introduced by Chen and Weth \cite{MR3995092}).
More recently, this was sharpened in~\cite{Lee2026, hauer2023characterization}, where the
boundedness assumption was removed at the expense of working within the framework of
Lizorkin distributions.

A common thread runs through all of these results: the underlying symbol vanishes only
to first order on the unit sphere $S^{d-1}$. In turn, the solution set is built entirely
from single-layer distributions. None of them, however, accommodate symbols with a
higher-order zero. The simplest equivalence already lying beyond their reach is
\begin{equation*}
	(1-(-\Delta)^\sigma)^2 u = 0 \quad\Longleftrightarrow\quad (1+\Delta)^2 u = 0,
	\qquad \text{in $\R^d$,}
\end{equation*}
whose symbol $(1-|\cdot|^{2\sigma})^2$ has a second-order zero on $S^{d-1}$, here the
solutions are genuine multi-layer distributions.
This raises a natural question: 
\begin{center}
\emph{Can such a rigidity result be recovered for symbols that
vanish to arbitrary finite order on the sphere, and if so, how does that order govern the
structure of the solution set?}
\end{center}

In this note we answer this question affirmatively. We remove the first-order
restriction entirely and consider symbols which vanish to finite order on the unit
sphere. We show that the order of vanishing precisely determines the multi-layer structure of
the solution set. As a by-product, we extend the framework
of~\cite{Lee2026, hauer2023characterization} beyond Bernstein functions of the Laplacian
and the logarithmic Laplacian, removing the boundedness assumption for a substantially
broader class of operators.

Before presenting the main result, we introduce some notation. 

\subsection{Notation and Conventions}

We use the following standard notation throughout:

\begin{itemize}
	\item $\mathcal{D}(\Omega)$: space of smooth, compactly supported functions on an open set $\Omega\subset \R^d$ or a manifold, \\
	$\mathcal{D}'(\Omega)$: its dual (distributions).
	\item $\mathcal{S}(\R^d)$: Schwartz space of rapidly decreasing smooth functions, \\
	$\mathcal{S}'(\R^d)$: tempered distributions.
	\item $\mathcal{F}$, $\mathcal{F}^{-1}$: Fourier and inverse Fourier transforms, normalized as
	\begin{equation*}
		\mathcal{F}(f)(\xi) = \hat{f}(\xi)=\int_{\R^d} e^{-ix\cdot\xi} f(x)\,dx, \qquad
		\mathcal{F}^{-1}(g)(x) = \frac{1}{(2\pi)^d} \int_{\R^d} e^{ix\cdot\xi} g(\xi)\,d\xi.
		\end{equation*}
	\item $S^{d-1}$: unit sphere in $\R^d$. 
\end{itemize}

All dualities are denoted by $\langle \cdot, \cdot \rangle$ with appropriate subscripts when needed. The Fourier transform is extended to distributions by duality. 

\subsection{Main Result}

The main problem considered in this note is the following:
\begin{equation}\label{eqn:intro}
	P(i\nabla)u = 0, \quad \text{in $\R^d$,}
\end{equation}

To properly define the notion of solution of~\eqref{eqn:intro} we need to introduce the Lizorkin space and distributions as well as the notion of admissible symbols. 

\begin{definition}\label{def:distributions}
		We define the Lizorkin space $\mathcal{Z}(\mathbb{R}^d)$ as the subspace of $\mathcal{S}(\mathbb{R}^d)$ (endowed with the same topology) given by
	\begin{equation*}
	\mathcal{Z}(\mathbb{R}^d)=\{\varphi\in \mathcal{S}(\mathbb{R}^d):\int_{\mathbb{R}^d}x^{\alpha}\varphi(x)\,dx=0, \text{ for all $\alpha \in \mathbb{N}^{d}_{0}$}\}.
	\end{equation*}
	Also, we denote 
	$\mathcal{Z}'(\mathbb{R}^d)$ as the set of continuous linear functionals on $\mathcal{Z}(\mathbb{R}^d)$. We call the space $\mathcal{Z}(\mathbb{R}^d)$ the Lizorkin space and $\mathcal{Z}'(\mathbb{R}^d)$ as the space of Lizorkin distributions. We also introduce
\begin{equation*}
		\Psi(\R^d):=\{\psi \in \mathcal{S}(\R^d):\partial^\alpha_{\xi} \psi(0)=0, \text{for all $\alpha \in \N^d_0$}\}.
	\end{equation*}
	We denote $\Psi'(\R^d)$ to be the space of linear continuous functionals on $\Psi(\R^d)$.
\end{definition}

\begin{definition}[Admissible symbols]
	Let $P:\R^d\rightarrow \mathbb{C}$ be a measurable function such that  
	\begin{itemize}
	\item[a)] $P$ is smooth on $\R^d\setminus \{0\}$, 
	\item[b)] there exists $l_{\gamma}\in \mathbb{N}$ and $C_{\gamma}>0$ such that 
			\begin{equation*}
				|\partial_{\xi}^\gamma P(\xi)|\leq C_{\gamma} 
				\begin{cases}
					|\xi|^{l_\gamma},& \quad \text{for $|\xi|\geq 1$},\\
					|\xi|^{-l_\gamma},& \quad \text{for $0<|\xi|<1$. }
				\end{cases}
			\end{equation*}
		\end{itemize}
	Then, we define the operator $P(i\nabla):\mathcal{Z}'(\R^d)\rightarrow \mathcal{Z}'(\R^d)$ as 
	\begin{equation*}
		\langle P(i\nabla )u, \varphi \rangle_{\mathcal{Z}'(\R^d), \mathcal{Z}(\R^d)}:=\langle u, P(-i\nabla )\varphi \rangle_{\mathcal{Z}'(\R^d), \mathcal{Z}(\R^d)},  \quad \text{for $(u,\varphi) \in \mathcal{Z}'(\R^d)\times \mathcal{Z}(\R^d)$.}
	\end{equation*}
	where
	\begin{equation*}
		P(-i\nabla)\varphi (x)=\frac{1}{(2\pi)^{d}}\int_{\R^d} e^{ix\cdot \xi}P(\xi)\hat{\varphi}(\xi)\,d\xi,\quad \text{for $x\in \R^d$. }
	\end{equation*}
	For the operator $P(i\nabla)$, we call $P$ the corresponding admissible symbol. 
\end{definition}

There are a few things to unpack in the context of the above definitions. 

Firstly, from~\cite[Lemma 2.4]{Kurokawa2002}, it is shown that the operator $P(i\nabla)$ is well-defined on the space of Lizorkin distributions. This is ultimately why we consider these distributions since they are compatible with the class of admissible symbols. To further emphasize this point, we utilize the logarithmic Laplacian $P(i\nabla)=\log(-\Delta)$, as a motivating example. Observe that if one was to naively define 
\begin{equation*}
	\langle\log(-\Delta)u,\varphi\rangle_{\mathcal{S}'(\R^d),\mathcal{S}(\R^d)}:=	\langle u,\log(-\Delta)\varphi\rangle_{\mathcal{S}'(\R^d),\mathcal{S}(\R^d)},\quad \text{for $(u,\varphi) \in \mathcal{S}'(\R^d)\times \mathcal{S}(\R^d)$},
\end{equation*}
one would run into an issue since 
\begin{equation*}
	\log(-\Delta)\varphi\notin \mathcal{S}(\R^d),
\end{equation*}
in general. It should be mentioned that polynomials are the problematic distributions for the logarithmic Laplacian. In turn, it requires for us to consider a notion of solution that allows us to work modulo polynomials. To be more precise, observe that  
	\begin{equation*}
	\langle u+ p,\varphi \rangle_{\mathcal{S}'(\R^d),\mathcal{S}(\R^d)}=        \langle u,\varphi \rangle_{\mathcal{S}'(\R^d),\mathcal{S}(\R^d)}, \quad \text{for all $(u,\varphi) \in \mathcal{S}'(\R^d)\times \mathcal{Z}(\R^d),$}
	\end{equation*}
	for all $p\in \mathcal{P}(\R^d)$ (space of polynomials). From this, we have that 
	$\mathcal{Z}'(\R^d)$ can be identified with the quotient space $\mathcal{S}'(\R^d)/\mathcal{P}(\R^d)$ and every element in $\mathcal{S}'(\R^d)$ can be identified as an element of $\mathcal{Z}'(\R^d)$ via projection. Every element in $\mathcal{Z}'(\R^d)$ can be uniquely identified (up to a polynomial) to an element in $\mathcal{S}'(\R^d)$. When we wish to emphasize that an element of $\mathcal{Z}'(\R^d)$ is an equivalence class we will use the notation $[u]$ as opposed to $u$. For further details on the Lizorkin space and distributions, we direct the reader to the book of Samko~\cite{MR1918790}. 

	We now introduce the notion of solution of~\eqref{eqn:intro}. 

\begin{definition}\label{def:notion_of_soln}
	We say that $u\in \mathcal{Z}'(\R^d)$ is a solution to~\eqref{eqn:intro} if 
	\begin{equation*}
		\langle u, P(-i\nabla)\varphi \rangle_{\mathcal{Z}'(\R^d),\mathcal{Z}(\R^d)}=0, 
	\end{equation*}
	for all $\varphi \in \mathcal{Z}(\R^d)$. We also identify $u=[u]\in \mathcal{Z}'(\R^d)$ with an element in $\mathcal{S}'(\R^d)$ such that $0\notin\supp(\hat{u})$. 
\end{definition}

\begin{remark}
	The last assumption can be viewed as a convenient choice for us to uniquely identify $\hat{u}\in \mathcal{S}'(\R^d)$ with a non-polynomial part. However, we need to assume this since it will become convenient for us to utilise the fact that $\R^d\setminus\{0\}$ and $(0,\infty)\times S^{d-1}$ are diffeomorphic. 
\end{remark}

To precisely state the main result, we need to introduce the space of single and multi layer distributions on $S^{d-1}$, cf.~\cite[Chapter 6]{MR2000535}. 

\begin{definition}
	\begin{itemize}
		\item We say that $T\in \mathcal{D}'(\R^d)$ is a \textbf{single layer distribution on the sphere}, if there exists a $\tilde{T}\in \mathcal{D}'(S^{d-1})$ such that 
	\begin{equation*}
		\langle T,\varphi\rangle_{\mathcal{D}'(\R^d),\mathcal{D}(\R^d)}=		\langle \tilde{T},\varphi|_{S^{d-1}}\rangle_{\mathcal{D}'(S^{d-1}),\mathcal{D}(S^{d-1})},
	\end{equation*}
	for all $\varphi \in \mathcal{D}(\R^d)$. We denote $\mathcal{SL}(S^{d-1})$ to be the set of single layer distributions on the sphere. 
	\item We say that $T\in \mathcal{D}'(\R^d)$ is a \textbf{multi-layer distribution (of order up to $N$) on the sphere} if $T$ can be expressed as:
	\begin{equation*}
		T=\sum_{k=0}^{N}\partial_r^k \tilde{T}_{k},
	\end{equation*}
	where $\tilde{T}_k$ is a single layer distribution on the sphere and $\partial_r$ is the radial derivative given by:
	\begin{equation*}
		\partial_r=\sum_{j=1}^{d}\frac{x_j}{|x|}\partial_{x_j}. 
	\end{equation*}
	We denote $\mathcal{ML}_N(S^{d-1})$ to be the set of multi-layer distributions (of order up to $N$). 
	\end{itemize}
\end{definition}

We are now ready to state the main result of this note. 

\begin{theorem}[Main result]\label{thm:main-intro}
Let $P$ be an admissible symbol in the sense of Definition~\ref{def:distributions},
and assume in addition that
\begin{enumerate}
  \item[(i)] the zero set of $P$ is exactly the unit sphere,
  \[
    \{\xi\in\R^d : P(\xi)=0\}=S^{d-1};
  \]
  \item[(ii)] $P$ vanishes to finite order on $S^{d-1}$: there is a largest
  integer $N\in\N_0$ such that
  \begin{equation}\label{eqn:flatness_condition}
    \lim_{r\to 1}\frac{P(r\omega)}{(r-1)^{N}}=0, \qquad\text{for all }\omega\in S^{d-1}.
  \end{equation}
\end{enumerate}
Then $u\in\mathcal{Z}'(\R^d)$ solves $P(i\nabla)u=0$ if and only if $\hat{u}\in\mathcal{ML}_N(S^{d-1})$.
\end{theorem}

\begin{remark}
An equivalent way of stating this result is that $u$ solves~\eqref{eqn:intro} if and only if $u$ solves $(1+\Delta)^{N+1} u=0$ in $\mathcal{S}'(\R^d)$, c.f.~\cite[Exemple, page 102]{MR209834}. Note we are able to say that $(1+\Delta)^N u=0$ in $\mathcal{S}'(\R^d)$, as opposed to $\mathcal{Z}'(\R^d)$, since we assume that $0\notin \supp(\hat{u})$. 
\end{remark}

\begin{remark}
	 In recent years there has been considerable interest in Liouville theorems for operators of L\'evy type cf.~\cite{MR4884561,MR4757484,MR4774765,MR4472712,MR4149690}. It is natural to consider whether the methods here can be applied to that setting. Unfortunately, there are some drawbacks since we typically work modul polynomials which are an important class of solutions to consider for Liouville theorems for L\'evy operators. Nevertheless, we should remark that ideas relating to division problems and generalized products have been considered in \cite{MR4884561,MR4774765}. 
\end{remark}

\begin{remark}
It is worthwhile to bring up some examples corresponding to our main result. It should be noted that our theory covers the m-th logarithmic Laplacian $\log(-\Delta)^m$ introduced by Huyuan Chen~\cite{MR5009321}. It is already noted in~\cite{Lee2026} that $\log(-\Delta):\mathcal{Z}'(\R^d)\rightarrow \mathcal{Z}'(\R^d)$ is well-defined. Hence, by induction, we have that the $\log(-\Delta)^m$ is also well-defined as well as satisfying a flatness condition near $S^{d-1}$ of order $m$. Other examples include $(f(-\Delta) -1)^m$ where $f$ is a non-constant Bernstein function.
\end{remark}

\subsection{Outline of the proof}

By the Fourier transform, it becomes apparent that in order to solve~\eqref{eqn:intro} it is equivalent to solve the following division problem:
\begin{equation}\label{eqn:division_problem}
	P(-\xi) \hat{u}=0, \quad \text{in $\Psi'(\R^d)$.}
\end{equation}
Namely, we need to classify all $\hat{u}\in {\Psi}'(\R^d)$ such that 
\begin{equation}\label{eqn:division_problem_distribution}
	\langle \hat{u},P\psi\rangle_{\Psi'(\R^d), \Psi(\R^d)}=0, \quad \text{for all $\psi \in \Psi(\R^d)$. }
\end{equation}

To show the forward direction --- that if $\hat{u}\in \mathcal{ML}_N(S^{d-1})$ then $\hat{u}$ solves~\eqref{eqn:division_problem} --- we proceed as follows:

\begin{enumerate}
	\item[Step 1.]
	We verify that the support of $\hat{u}$ must be contained in $S^{d-1}$. From this, we can replace the condition $\psi \in \Psi(\R^d)$, in \eqref{eqn:division_problem_distribution}, to $\psi \in \mathcal{D}(\R^d)$ such that $\supp(\psi)\subset \{1-\epsilon<|\xi|<1+\epsilon\}$ for some sufficiently small $\epsilon>0$.  
	\item[Step 2.]
	We show that $\partial_{r}^jP|_{S^{d-1}}=0$, for $j=0,..,N$. This will assist us in proving the first direction of Theorem~\ref{thm:main-intro} which we do in Step 3.  
	\item[Step 3.]
	We show that if $\hat{u}\in \mathcal{ML}_N(S^{d-1})$, then we have that $\hat{u}$ is a solution of~\eqref{eqn:division_problem}.
\end{enumerate}	
	
We now consider the other direction. Before we explicitly outline our steps, note the following:

If we fix $\tau \in \mathcal{D}(S^{d-1})$ and define $\hat{u}_{\text{Sph}}^\tau\in \mathcal{D}'((0,\infty))$:

\begin{equation}\label{eqn:spherical_coordinates}
	\langle \hat{u}_{\text{Sph}}^\tau, \rho\rangle_{\mathcal{D}'((0,\infty)), \mathcal{D}((0,\infty))}:=\langle \hat{u},\rho \tau\rangle_{\mathcal{D}'(\R^d),\mathcal{D}(\R^d)}, \text{for $\rho \in \mathcal{D}((0,\infty))$}, 
\end{equation}
where 
\begin{equation*}
	(\rho\tau)(\xi):=\rho(|\xi|)\tau\left ( \tfrac{\xi}{|\xi|}\right ), \quad \text{ for $\xi \in \R^d$,}
\end{equation*}
then from Step 1, we have that $\supp(\hat{u}_{\text{Sph}}^\tau)=\{1\}$. Necessarily, by the structure theorem of distributions supported at a point~\cite[Theorem 2.3.4]{MR1996773}, we have that
\begin{equation}\label{eqn:point_structure}
	\langle \hat{u}_{\text{Sph}}^{\tau}, \rho \rangle_{\mathcal{D}'((0,\infty)), \mathcal{D}((0,\infty))}=\sum_{k=0}^{M} c_k(\tau)\partial_{r}^{k}\rho(1),
\end{equation}
for some $M\in \N_0$. Note that if we choose $\rho$ such that $\partial_{r}^j\rho(1)=0$ for $j\ne k$ then one can see from the definition given in~\eqref{eqn:spherical_coordinates} that each coefficient $\tau\mapsto c_k(\tau)$ defines a distribution $c_k\in \mathcal{D}'(S^{d-1})$.

In turn, ~\eqref{eqn:point_structure} already exhibits $\hat{u}$, in spherical coordinates, as a multi-layer distribution of order at most $M$; a priori, this order may exceed $N$. To conclude that $\hat{u}\in \mathcal{ML}_N(S^{d-1})$, it therefore remains to show that the higher-order coefficients $c_k$ vanish for every $k>N$. We carry this out in the following two steps.

\begin{enumerate}
	\item[Step 4.]
	We show that the action of $\hat{u}$ on a test function is determined by its radial Taylor jet of order $N$ along $S^{d-1}$. Concretely, we prove the following.
		
	Fix $\tau_{k}\in \mathcal{D}(S^{d-1})$ for $k=0,..,N$ and let $\varphi \in \mathcal{D}(\R^d)$ such that 
	 \begin{equation*}
		\partial_{r}^k\varphi|_{S^{d-1}}=\tau_{k}\in \mathcal{D}(S^{d-1}), \quad \text{for $k=0,..,N$.}
	 \end{equation*} 
	By taking $r=|\xi|$ and $\omega = \tfrac{\xi}{|\xi|}$, denote $E\tau:\R^d\rightarrow \mathbb{C}$ by
	\begin{equation*}
		(E\tau)(r\omega) := \chi(r)\sum_{k=0}^{N} \frac{(r-1)^k}{k!}\tau_k(\omega), \quad \text{for $r>0$ and $\omega\in S^{d-1}$,}
	\end{equation*}
	where $\chi$ is a non-negative smooth bump function with support contained in $\{1-\epsilon<r<1+\epsilon\}$ such that $\chi \equiv 1$ on $\{1-\tfrac{\epsilon}{2}<r<1+\tfrac{\epsilon}{2}\}$. Then we have that
	 \begin{equation*}
			\langle\hat{u},\varphi \rangle_{\mathcal{D}'(\R^d), \mathcal{D}(\R^d)}=\langle\hat{u},E\tau \rangle_{\mathcal{D}'(\R^d), \mathcal{D}(\R^d)}. 
	 \end{equation*} 
	 \item[Step 5.]
	By Step 4, the pairing $\langle\hat{u}_{\text{Sph}}^\tau,\rho\rangle$ depends only on the radial jet $\partial_r^k\rho(1)$ for $k=0,\dots,N$. Comparing this with the representation~\eqref{eqn:point_structure} forces $c_k=0$ for every $k>N$, and hence $\hat{u}=\sum_{k=0}^{N}\partial_{r}^{k}\tilde{T}_{k}\in \mathcal{ML}_N(S^{d-1})$. Together with Step 3, this completes the proof of Theorem~\ref{thm:main-intro}.
\end{enumerate}

\begin{remark}
It is worthwhile to remark that the proofs considered in~\cite{MR4709301, MR4546886, MR4845366, hauer2023characterization, MR3477075, Lee2026} are genuinely multi-dimensional. Here we are somewhat more consistent with the approach considered by Fall and Weth~\cite{MR3511811} in the sense that we utilize the structure theorem of distributions supported at a point.  
\end{remark}

\section{Proof of Theorem~\ref{thm:main-intro}}

We follow the steps as outlined in the previous section. 

\subsection*{Step 1.}

We are required to show that, for all $\varphi \in \mathcal{D}(\R^d)$ such that $S^{d-1}\cap \supp(\varphi) =\emptyset$, we have that 
\begin{equation*}
\langle \hat{u},\varphi \rangle_{\Psi'(\R^d),\Psi(\R^d)}=0. 
\end{equation*}

Since we assume that $0\notin \supp(\hat{u})$, we can identify $\hat{u}$ as an element of $\mathcal{D}'(\R^d)$. 

The function $\omega:\R^d\rightarrow \mathbb{C}$ given by 
\begin{equation*}
	\omega(\xi)=\begin{cases}
		\frac{1}{P(\xi)}\varphi(\xi), \quad &\text{for $\xi \in \supp(\varphi)$}, \\
		0, \quad &\text{otherwise,}
	\end{cases} 
\end{equation*}
belongs to $\Psi(\R^d)$ since $P$ is smooth with non-zero derivatives (on $\supp(\varphi)$). In turn, we have that 
\begin{equation*}
	\begin{split}
	\langle \hat{u},\varphi \rangle_{\Psi'(\R^d),\Psi(\R^d)}
	&=\langle \hat{u},\varphi \rangle_{\mathcal{D}'(\R^d),\mathcal{D}(\R^d)},\\
	&=\langle \hat{u},P\omega \rangle_{\mathcal{D}'(\R^d),\mathcal{D}(\R^d)},\\
	&=\langle \hat{u},P\omega \rangle_{\Psi'(\R^d),\Psi(\R^d)},\\
	&=0.
	\end{split}
\end{equation*}
In the last line, we use~\eqref{eqn:division_problem_distribution}. 
In turn, we have verified Step 1. 

\subsection*{Step 2.}

Denote $f(r):=P(r\omega)$, for $r>0$ and fixed $\omega \in S^{d-1}$. Since $f$ is smooth, we have, by the Taylor remainder theorem:
\begin{equation*}
f(r)=\sum_{j=0}^{N}\frac{f^{(j)}(1)}{j!}(r-1)^j + O(|r-1|^{N+1}), \quad \text{for $r$ near $1$. }
\end{equation*}
As a consequence of~\eqref{eqn:flatness_condition}, it becomes necessary that $f^{(j)}(1)=0$, for $j=0,...,N$, which completes Step 2. 

\subsection*{Step 3.}

We assume that $\hat{u}\in \mathcal{ML}_{N}(S^{d-1})$ of the form:

\begin{equation*}
	\hat{u}=\sum_{k=0}^{N}\partial_r^k\tilde{T}_{k}, 
\end{equation*}
where $\Tilde{T}_{k} \in \mathcal{SL}(S^{d-1})$. By linearity, it is sufficient to show that 
\begin{equation*}
	\langle \partial^{k}_{r}\Tilde{F},P\psi\rangle_{\Psi'(\R^d), \Psi(\R^d)}=0, \quad \text{for $k=0,..,N$ and $\psi \in \Psi(\R^d)$,}
\end{equation*}
where $\Tilde{F}\in \mathcal{SL}(S^{d-1})$, denote $F\in \mathcal{D}'(S^{d-1})$ to be the distribution on the sphere corresponding to $\Tilde{F}\in \mathcal{SL}(S^{d-1})$. From Step 1, it is sufficient to take $\psi \in \mathcal{D}(\R^d)$ such that $\supp(\psi)\subset \{1/2<|\xi|<3/2\}$. By direct computation, we have that 
\begin{equation*}
	\begin{split}
	\langle \partial^{k}_{r}\Tilde{F},P\psi\rangle_{\Psi'(\R^d), \Psi(\R^d)}
	&=\langle \Tilde{F},(\partial_{r}^{*})^{k}(P\psi)\rangle_{\Psi'(\R^d), \Psi(\R^d)}, \\
	&=\langle F,(\partial_{r}^{*})^{k}(P\psi)|_{S^{d-1}}\rangle_{\mathcal{D}'(S^{d-1}), \mathcal{D}(S^{d-1})}, 
	\end{split}
\end{equation*}
where $\partial_r^*$ denotes the adjoint of $\partial_r$. To complete Step 3, we verify that $(\partial_{r}^{*})^{k}(P\psi)|_{S^{d-1}}=0$. 

Note that $(\partial_{r}^{*})^{k}f=(-1)^kr^{-(d-1)}\partial^k_{r}(r^{d-1}f)$. Hence, by an application of Step 2 and the product rule, we have that $(\partial_{r}^{*})^{k}(P\psi)|_{S^{d-1}}=0$. 
\subsection*{Step 4.} 

To verify the statement, observe the following decomposition:
\begin{equation*}
	\begin{split}
		\varphi(r\omega)-(E\tau)(r\omega)=\chi(r)\varphi(r\omega)-(E\tau)(r\omega)+(1-\chi(r))\varphi(r\omega),\quad \text{for $r>0$ and $\omega \in S^{d-1}$. }
	\end{split}
\end{equation*}
 The last term will be annihilated by $\hat{u}$ since $\supp(\hat{u})\subset S^{d-1}$. In turn, we focus our attention on $\chi(r)\varphi(r\omega)-(E\tau)(r\omega)$. By the Taylor expansion, we have that:

\begin{equation*}	
\begin{split}
	\chi(r)\varphi(r\omega)-(E\tau)(r\omega)
&=
(r-1)^{N+1} \left (\frac{\chi(r)}{N!}
\int_{0}^{1}
(1-t)^N\,
\partial_r^{N+1}
\varphi\left((1+t(r-1))\omega\right)\,dt \right ), \\
&=
P(r\omega)\underbrace{\left (\frac{(r-1)^{N+1}}{P(r\omega)} \right )}_{\tfrac{1}{Q(r,\omega)}}\underbrace{\left (\frac{\chi(r)}{N!}
\int_{0}^{1}
(1-t)^N\,
\partial_r^{N+1}
\varphi\left((1+t(r-1))\omega\right)\,dt \right )}_{R(r,\omega)}. 
\end{split}
\end{equation*}
Note that $R$ is clearly smooth and compactly supported. Under the assumption that $1/Q$ is smooth on the support of $\chi$ then we have that 
\begin{equation*}
	\langle \hat{u}, \varphi -E\tau\rangle_{\mathcal{D}'(\R^d),\mathcal{D}(\R^d)}=\langle \hat{u}, P\left (\tfrac{1}{Q}\right )R\rangle_{\mathcal{D}'(\R^d),\mathcal{D}(\R^d)}=0, 
\end{equation*}
from~\eqref{eqn:division_problem_distribution}. To show that $1/Q$ is smooth (on a sufficiently small enough neighbourhood of $S^{d-1}$), we need to show that $Q(r,\omega):=\frac{P(r\omega)}{(r-1)^{N+1}}$, is smooth and is non-zero (on some neighbourhood of $S^{d-1}$). 

From~\eqref{eqn:flatness_condition}, Step 2 and the Taylor expansion of $P$, we have that 
\begin{equation*}
	Q(r,\omega)= \frac{1}{N!}\int_0^{1}(1-t)^N\partial_{r}^{N+1}P((1+t(r-1))\omega)\,dt, \quad \text{for $r>0,\omega \in S^{d-1}$. }
\end{equation*}

From this, we have that $Q$ is smooth. Moreover, by continuity of $Q$ and the definition of $N$, we have that $Q$ is non-zero for $r\in (1-\epsilon, 1+\epsilon)$ for a sufficiently small enough $\epsilon>0$. Hence, by the quotient rule, we have that $1/Q$ is smooth. This conclude Step 4. 

\subsection*{Step 5.}

Recall from~\eqref{eqn:spherical_coordinates}--\eqref{eqn:point_structure} that, for every $\tau\in \mathcal{D}(S^{d-1})$, the distribution $\hat{u}_{\text{Sph}}^{\tau}$ is supported at $\{1\}$ and satisfies
\begin{equation*}
	\langle \hat{u}_{\text{Sph}}^{\tau}, \rho \rangle=\sum_{k=0}^{M} c_k(\tau)\,\partial_{r}^{k}\rho(1), \qquad c_k\in \mathcal{D}'(S^{d-1}),
\end{equation*}
for some $M\ge N$. It thus suffices to show that $c_k=0$ for every $k>N$.

Fix $j$ with $N<j\le M$ and $\tau\in \mathcal{D}(S^{d-1})$. By Step 4, on $\{1-\epsilon<|\xi|<1+\epsilon\}$ we have the factorization $P(r\omega)=(r-1)^{N+1}Q(r,\omega)$ with $Q$ smooth and non-vanishing. Let $\rho\in \mathcal{D}((0,\infty))$ be such that $\partial_{r}^{k}\rho(1)=0$ for $k=0,..,N$, and set
\begin{equation*}
	\psi(r\omega):=\frac{\rho(r)\,\tau(\omega)}{P(r\omega)}=\frac{\rho(r)}{(r-1)^{N+1}Q(r,\omega)}\,\tau(\omega).
\end{equation*}
Then $\psi\in \mathcal{D}(\R^d)$ with $\supp(\psi)\subset \{1-\epsilon<|\xi|<1+\epsilon\}$, since $\rho(r)/(r-1)^{N+1}$ and $1/Q$ are smooth.

Since $P\psi=\rho\tau$, we have, from~\eqref{eqn:division_problem_distribution} and~\eqref{eqn:point_structure},
\begin{equation*}
	0=\langle \hat{u},P\psi\rangle=\langle \hat{u}_{\text{Sph}}^{\tau},\rho\rangle=\sum_{k=N+1}^{M}c_k(\tau)\,\partial_r^k\rho(1),
\end{equation*}
the sum starting at $k=N+1$ because $\partial_r^k\rho(1)=0$ for $k=0,.., N$. Choosing $\rho$ appropriately gives $c_j(\tau)=0$.

Since $j>N$ and $\tau$ were arbitrary, $c_k=0$ for all $k>N$, so $\hat{u}_{\text{Sph}}^{\tau}=\sum_{k=0}^{N}c_k(\tau)\partial_r^k\delta_1$. Equivalently, $\hat{u}=\sum_{k=0}^{N}\partial_r^k\tilde{T}_k\in \mathcal{ML}_N(S^{d-1})$, where $\tilde{T}_k\in \mathcal{SL}(S^{d-1})$ is the single layer corresponding to $c_k$. Hence, utilising~\cite[Ch.~IV, Thm.~III]{MR209834} we have our result since we test against products of $\mathcal{D}((0,\infty))$ and $\mathcal{D}(S^{d-1})$ which are dense in $\mathcal{D}(\R^{d}\setminus\{0\})$. This completes the proof of Theorem~\ref{thm:main-intro}. 

\begin{remark}
Note that while the result for~\cite[Ch.~IV, Thm.~III]{MR209834} isn't stated for manifolds it is clarified on page 108 that without modification this can be done. 
\end{remark}

\section*{Acknowledgements}

The author is very appreciative of Clio Cresswell and thanks her for her encouragement and advice.

\end{document}